
\input eplain.tex

\vsize=21cm 
\voffset=1.8cm
\baselineskip=11.5pt 

\def\rightheadline{\hfil{\small A note on the mean value of
the zeta and $\scriptstyle L$-functions. XV}
\hfil\tenrm\folio}
\def\leftheadline{\tenrm\folio\hfil{\scsc Y. Motohashi}\hfil}
\def\emptyheadline{\hfil}
\headline{\ifnum\pageno=1 \emptyheadline\else
\ifodd\pageno \rightheadline \else \leftheadline\fi\fi}

\def\firstpage{\hss{\vbox to 2cm{\vfil\hbox{\rm\folio}}}\hss}
\def\emptyfootline{\hfil}
\footline{\ifnum\pageno=1\firstpage\else
\emptyfootline\fi}

\input amssym.def
\input amssym.tex

\font\small=cmr8
\font\csc=cmcsc10
\font\scsc=cmcsc10 at 8pt
\font\title=cmbx10
\font\stitle=cmmib10

\font\csc=cmcsc10
\font\title=cmbx10 at 12pt
\font\stitle=cmmib10 at 12pt

\font\teneusm=eusm10
\font\seveneusm=eusm7
\font\fiveeusm=eusm5
\newfam\eusmfam
\def\eusm{\fam\eusmfam\teneusm}
\textfont\eusmfam=\teneusm
\scriptfont\eusmfam=\seveneusm
\scriptscriptfont\eusmfam=\fiveeusm

\def\varGamma{{\mit \Gamma}}
\def\Re{{\rm Re}\,}

\def\txt#1{{\textstyle{#1}}}
\def\scr#1{{\scriptstyle{#1}}}
\def\r#1{{\rm #1}}
\def\B#1{{\Bbb #1}}
\def\e#1{{\eusm #1}}

\def\sgn{{\rm sgn}}

\singlecolumn
\centerline{\title A note on the mean value of the zeta and 
{\stitle L\/}-functions. XV} 
\vskip 1cm
\centerline{By Yoichi {\csc Motohashi}${}^{\ast)}$}
\vskip 1cm 
\hsize=16.5cm
{\bf Abstract:} The aim of the present article is to render   
the spectral theory of mean values of automorphic $L$-functions --
in a unified fashion.
This is an outcome of the investigation commenced with the
parts XII and XIV, where a framework was laid on the
basis of the theory of automorphic representations and a general
approach to the mean values was envisaged.  We restrict ourselves to
the situation offered by the full modular group, solely for the sake
of simplicity. Details and
extensions are to be published elsewhere.
\par
\hoffset=-0.2cm
\doublecolumns
\hsize=227.38682pt
\medskip
\footnote{}{\small 2000 Mathematics Subject 
Classification:  11F70
\par
${}^{\ast)}$ Department of Mathematics, 
Nihon University, Surugadai,
Tokyo 101-8308.} 
{\bf 1.} To begin with, we stress that all notations and conventions
will stay effective once introduced. 
\smallskip
Let $\r{G}=\r{PSL}(2,{\Bbb
R})$ and $\varGamma=\r{PSL}(2,{\Bbb Z})$. Write
$$
\eqalign{
&\r{n}[x]=\left[\matrix{1&x\cr&1}\right],\quad
\r{a}[y]=\left[\matrix{\sqrt{y}&\cr&1/\sqrt{y}}\right],\cr
&\r{k}[\theta]=\left[\matrix{\phantom{-}\cos\theta&
\sin\theta\cr-\sin\theta &\cos\theta}\right],
}
$$
and $\r{N}=\left\{\r{n}[x]: x\in\B{R}\right\}$,
$\r{A}=\left\{\r{a}[y]: y>0\right\}$,
$\r{K}=\left\{\r{k}[\theta]:\theta\in\B{R}/\pi\B{Z}\right\},
$ 
so that $\r{G}=\r{NAK}$ is the Iwasawa decomposition 
of the Lie group $\r{G}$.  We read it as 
$\r{G}\ni\r{g}=\r{n}\r{a}\r{k}=\r{n}[x]\r{a}[y]\r{k}[\theta]$.
The coordinate $(x,y,\theta)$ will retain this definition. 
The center of the universal enveloping algebra 
of $\r{G}$ is the polynomial ring on the Casimir element  
$\Omega=y^2\!\left(\partial_x^2+
\partial_y^2\right)-y\partial_x\partial_\theta$.
The Haar
measures on the groups
$\r{N}$, $\r{A}$, $\r{K}$, $\r{G}$ are normalized, respectively, by
$d\r{n}=dx$,
$d\r{a}=dy/y$, $d\r{k}=d\theta/\pi$,
$d\r{g}=d\r{n}d\r{a}d\r{k}/y$,  with Lebesgue measures $dx$,
$dy$, $d\theta$. 
\par
The space $L^2(\varGamma\backslash\r{G})$ is composed of
all left $\varGamma$-automorphic functions on $\r{G}$, 
vectors for short, which are square integrable over
$\varGamma\backslash\r{G}$ against
$d\r{g}$. Elements of $\r{G}$ act unitarily on vectors from the right.
We have the orthogonal decomposition
$$ 
L^2(\varGamma\backslash\r{G})=\B{C}\!\cdot\!1\txt{\bigoplus}
{}^0\!L^2(\varGamma\backslash\r{G})\txt\bigoplus
{}^e\!L^2(\varGamma\backslash\r{G})
$$ 
into invariant subspaces. Here ${}^0\!L^2$ is the cuspidal subspace
spanned by vectors whose Fourier expansions with respect to the left
action of $\r{N}$ have vanishing constant terms. The subspace
${}^e\!L^2$ is spanned by integrals of Eisenstein series. 
Invariant subspaces of $L^2(\varGamma\backslash\r{G})$ and
$\varGamma$-automorphic representations of
$\r{G}$ are interchangeable  concepts, and 
we shall refer to them in a
mixed way. 
\par
The cuspidal subspace decomposes into irredu\-cible subspaces
$$ 
{}^0\!L^2(\varGamma\backslash\r{G})=\overline{\txt{\bigoplus} V}.
$$ 
The operator $\Omega$ becomes a constant multiplication in each
$V$ so that 
$\Omega|_{V^\infty}=\!\left(\nu^2_V-{1\over4}\right)\cdot1$, 
where $V^\infty$ is the set of all infinitely differentiable 
vectors in $V$. Under our present supposition, $V$ belongs to either
the unitary principal series or the discrete series; accordingly,
we have $\nu_V\in i\B{R}$ or $\nu_V=\ell-{1\over2}$
$(1\le\ell\in\B{Z})$. 
\par
The right action of $K$ induces the decomposition of
$V$ into $K$-irreducible subspaces
$$ 
V=\overline{{\mathop\txt{\bigoplus}_{p=-\infty}^\infty} V_p}\,,\quad
\dim V_p\le 1.
$$ 
If it is not trivial, $V_p$ is spanned by a
$\varGamma$-automorphic function on which the right translation by
$\r{k}[\theta]$ becomes the multiplication by the factor
$\exp(2ip\theta)$. It is called a $\varGamma$-automorphic form of
spectral parameter $\nu_V$ and weight $2p$.
\par
Let $V$ be in the unitary principal series.  Then $\dim V_p=1$ for all
$p\in\B{Z}$ and there exists a complete orthonormal system
$\{\varphi_p\in V_p:\,p\in\B{Z}\}$ of $V$ such that
$$
\varphi_p(\r{g},V)=\sum_{\scr{n=-\infty}\atop\scr{n\ne0}}
^\infty{\varrho_V(n)\over\sqrt{|n|}}
\e{A}^{\sgn(n)}\phi_p(\r{a}[|n|]\r{g};\nu_V),
$$ 
where $\phi_p(\r{g};\nu)=y^{\nu+{1\over2}}\exp(2ip\theta)$, and
$$
\e{A}^\delta\phi_p(\r{g};\nu)
=\int_{-\infty}^\infty
\exp(-2\pi i\delta x)\phi_p (\r{w}\r{n}[x]\r{g};\nu)dx,
$$ 
with $\r{w}=\r{k}\!\left[{1\over2}\pi\right]$ the Weyl element.
Our normalization is such that the coefficients
$\varrho_V(n)$ do not depend on the weight. Also
we may impose the Hecke invariance; in particular
$\varrho_V(-n)=\epsilon_V\varrho_V(n)$ with $\epsilon_V=\pm1$.
\par
We have skipped the discrete series. In the sequel
as well, we shall argue as if there were no discrete series
representations. This should not cause any confusion, as the
discussions below, especially those pertaining estimations, extend
readily to the discrete series. Nevertheless, it should be remarked
that the definition of $\varphi_p(\r{g},V)$ 
given above has to be modified for those $V$ in the discrete series
according to the normalization [2, $(2.21)$ and $(2.25)$].
\smallskip
Our discussion depends much on the uniform bounds for
$\e{A}\phi_p(\r{a}[y],\nu)$, $\e{A}=\e{A}^+$, such as [2, $(4.3)$ and
$(4.5)$]. In order to make our argument applicable to any cuspidal
representations, we derive from the latter a bound that is somewhat
weaker than the former but is still sufficient for our purpose; in
fact the proof of [2, $(4.5)$] works for all cases. We thus put
$$
\Gamma_p(s,\nu)=\int_0^\infty y^{s-1}\e{A}\phi_p(\r{a}[y],\nu)
{dy\over y},\quad \Re s>1;
$$
we have shifted the argument in [2, $(4.10)$] by $-{1\over2}$. Let us
assume that $\nu\in i\B{R}$. We divide the integral at $y
=|p|+|\nu|+1$. To the part with smaller $y$ we apply the fact that
$\e{A}^{\sgn(u)}\phi_p(\r{a}[|u|];\nu)$ is a unit vector in
$L^2(\B{R}\!^\times\!, d^\times\!/\pi\!)$, $d^\times\!u=du/|u|$ 
(see e.g., [2, Lemma 4]). Hence this part is
$\ll (|p|+|\nu|+1)^{\Re s-1}$. On the other  hand, by [2,
$(4.5)$] the remaining part is $\ll(|p|+|\nu|+1)^{\Re s-{1\over2}}$.
We then invoke the identity
$$
\Gamma_p(s,\nu)=4\pi\cdot
{\pi\Gamma_p(s+2,\nu)-p\Gamma_p(s+1,\nu)\over
(s-{1\over2})^2-\nu^2},
$$ 
which can be proved by applying the operator
$\e{D}_\nu=(d/dy)^2-(2\pi)^2-\left(\nu^2-
\txt{1\over4}\right)\!y^{-2}$ either to $y^s$ or to
$\e{A}\phi_p(\r{a}[y],\nu)$ in the last integral, on noting 
that the latter is a constant
multiple of the Whittaker function $W_{p,\nu}(4\pi y)$ (see [2,
(2.16)]). Then, by Mellin's inversion, we conclude that
$\e{A}\phi_p(\r{a}[y],\nu)
\ll y^{{1\over2}-\varepsilon}(|p|+|\nu|+1)^{2+\varepsilon}$
for any small $\varepsilon>0$. 
\par
The above works with $\nu=\ell-{1\over2}$, the discrete series, 
as well. We remark that $\e{D}_\nu$ is related to
the operator $\partial_\theta$ via the Kirillov map.
Hereafter we shall term  [2, $(4.5)$] and the bound thus obtained
the basic bounds.
\medskip
{\bf 2.} We now define the automorphic $L$-function
associated with an irreducible $\varGamma$-automorphic representation
$V$ by
$$
L_V(s)=\sum_{n=1}^\infty\varrho_V(n)n^{-s}.
$$
This converges absolutely for $\Re s>1$ and continues to an entire
function  which is of a polynomial order both in $s$ and in $\nu_V$ if
$\Re s$ is bounded. 
\par
We fix an $A$ among $V$'s, and consider the mean square
$$
\e{M}(A,g)=\int_{-\infty}^\infty\left|L\!\left(\txt{1\over2}+it\right)
\right|^2g(t)dt,\quad L=L_A,
$$
where the weight function
$g$ is assumed, for the sake of simplicity but without much loss
of generality, to be   even, entire, real on $\Bbb R$, and of fast
decay in any fixed horizontal strip. Our aim is 
to establish a full spectral decomposition of $\e{M}(A,g)$. 
Our method is applicable to any representation, though we shall deal
with only those $A$ in the unitary principal series. With a
minor modification, the fourth moment of the Riemann zeta-function
could be discussed equally. 
\par
We start as usual with the integral
$$
I(u,v;g)=\int_{-\infty}^\infty
\overline{L(\bar{u}+it)}L(v+it)g(t)dt.
$$
This is an entire function of $u,\,v$; in particular
$$
\e{M}(A,g)=I\!\left(\txt{1\over2},\txt{1\over2};g\right).
$$ 
In the region of absolute convergence we have
$$
I(u,v;g)={R(u+v)\over\zeta(2(u+v))}\hat{g}(0)+J(u,v;g)
+\overline{J(\bar{v},\bar{u};g)},
$$
where $R$ is the Rankin zeta-function attached to $A$, 
$\hat{g}$ the Fourier transform of $g$, and
$$
\eqalign{
&J(u,v;g)=\sum_{f=1}^\infty\sum_{n=1}^\infty
{\overline{\varrho(n)}\varrho(n+f)\over (2n+f)^{u+v}}\cr
&\cdot\left({\sqrt{n(n+f)}\over 2n+f}\right)^{2\alpha}
\tilde{g}(f/(2n+f);u,v), 
}
$$
where $\varrho=\varrho_A$ and
$$
\tilde{g}(x;u,v)
=2^{u+v+2\alpha}{\hat{g}(\log( (1+x)/(1-x)))\over (1-x)^{u+\alpha}
(1+x)^{v+\alpha}},
$$
with $0\le x\le1$.
Here $\alpha$ is a sufficiently large positive integer,
which is implicit throughout the sequel. Note that we
have  slightly modified our procedure given in [6].
\par
Let $g^*$ be the Mellin transform
of $\tilde{g}$. It is immediate to see that
$g^*(s;u,v)$ is of rapid decay with respect to
$s$, provided $\Re s$ and $u,v$ are bounded; moreover,
$g^*(s;u,v)/\Gamma(s)$ is entire over $\B{C}^3$.  Thus, by Mellin's
inversion, 
$$
\eqalign{
&J(u,v;g)\cr
&={1\over2\pi i}\int_{(\eta)}\!\left\{\sum_{f=1}^\infty f^{-s}
D_f(u+v-s)\right\}
g^*(s;u,v)ds,
}
$$
where $(\eta)$ is the line $\Re s=\eta>0$, and
$$
D_f(s)=\sum_{n=1}^\infty
{\overline{\varrho(n)}\varrho(n+f)\over (2n+f)^s}
\!\left({\sqrt{n(n+f)}\over 2n+f}\right)^{2\alpha}.
$$
Checking the convergence, we see that the condition
$\Re(u+v)>\max\{2,1+\eta\}$ is required here.
\medskip
{\bf 3.} Now, we need to have a full spectral decomposition of
$D_f(s)$ which is to yield a continuation, with a polynomial
growth, to the domain $\Re s>\varpi$ with a $\varpi<1$  so that 
$J(u,v;g)$  can also be continued to a neighbourhood of the 
point $\left({1\over2},{1\over2}\right)$. 
Recently, V. Blomer and G. Harcos [1] succeeded in establishing such
an assertion on $D_f(s)$. They
started with an old idea of ours, published in [6], to employ the
Kirillov model to pick up a favourable automorphic form in dealing
with the problem of the spectral decomposition of
$D_f(s)$, and proceeded one step further by a use of the theory of the
Sobolev norms with which they could derive the crucial polynomial
growth that we  had left open in [6]. 
\smallskip
In what follows, we shall take an alternative way by modifying
[6] with certain simple devices from [2], [5], and by reworking
our relevant unpublished notes. The present article can be read
independently of [1].
\par
Thus, let $A$, $\alpha$ be as above, and $\tau$ a
parameter in the right half plane; all implicit constants in the
sequel may depend on $A$, $\alpha$ and $\Re\tau$ at most. We apply
the inverse Kirillov map to the function $w(y,\tau)$ which is equal to
$y^{\alpha+{1\over2}}\exp(-\tau y)$ for
$y>0$, and vanishes for $y\le0$.  By [2, Lemma 4] or [6, Lemma
1], there exists an automorphic form $\Phi(\r{g},\tau)$ in 
the subspace $A$ such that
$$
\Phi(\r{n}[x]\r{a}[y],\tau)=
\sum_{n=1}^\infty{\varrho(n)\over\sqrt{n}}w(ny,\tau)
\exp(2\pi inx).
$$
More precisely,
$$
\Phi(\r{g},\tau)=\sum_{p=-\infty}^\infty a_p(\tau)
\varphi_p(\r{g}, A),
$$
where
$$
a_p(\tau)={1\over\pi}\int_0^\infty w_\alpha(y,\tau)
\overline{\e{A}\phi_p(\r{a}[y];\nu_A)}\,{dy\over y}.
$$
The function $\Phi(\r{g},\tau)$ is regular for $\Re \tau>0$. 
In fact, we have
$$
a_p(\tau)\ll (|\tau|+1)^{2\alpha}(|p|+1)^{-\alpha};
$$
and for any $\r{g}$ it holds that $\varphi_p(\r{g},A)\ll 
(|p|+1)^2$, by [2, $(4.5)$] applied to
$\e{A}\phi_p(\r{a}[y];\nu_A)$. 
\par
To prove the bound for $a_p(\tau)$, we use the operator $\e{D}_{\nu}$
again:  we have, for $p\ne0$,
$$
a_p(\tau)=-{1\over4\pi^2p}\int_0^\infty w_\alpha(y,\tau)
\overline{\e{D}_{\nu_A}\e{A}\phi_p(\r{a}[y];\nu_A)}\,dy,
$$
Integrate in parts, and repeat the procedure as long as
the resulting integrand vanishes at $y=0$, that is, $\alpha$ times
at least, and we get the bound. The case $p=0$ can
be skipped obviously. 
\medskip
{\bf 4.} Next, we put $\Psi(\r{g},\tau)=\Phi(\r{g},\tau)
\overline{\Phi(\r{g},\overline{\tau})}$; this is regular for
$\Re\tau>0$. We have, for any integer $f\ge0$,
$$
\eqalign{
&\int_0^1 \Psi(\r{n}[x]\r{a}[y],\tau)\exp(-2\pi ifx)dx\cr
&=y^{2\alpha+1}\sum_{n=1}^\infty\overline{\varrho(n)}
\varrho(n+f)(n(n+f))^\alpha\cr
&\hskip 1cm\cdot\exp(-(2n+f)\tau y).
}
$$
The Parseval formula implies that
the left side is equal to
$$
\eqalign{
&\sum_V{\varrho_V(f)\over\sqrt{f}}\Delta(fy,\tau;V)
\cr
+&\int_{(0)}
{f^{-\nu}\sigma_{2\nu}(f)\over
\sqrt{f}\zeta(1+2\nu)}\Delta(fy,\tau;\nu){d\nu\over4\pi i},
}
$$
where $V$ runs over all irreducible cuspidal representations, and 
$$
\eqalign{
\Delta(y,\tau;V)&=\sum_{p=-\infty}^\infty\langle\Psi(\cdot,\tau),
\varphi_p(\cdot,V)\rangle\e{A}\phi_p(\r{a}[y];\nu_V),\cr
\Delta(y,\tau;\nu)&=\sum_{p=-\infty}^\infty\langle\Psi(\cdot,\tau),
E_p(\cdot,\nu)\rangle\e{A}\phi_p(\r{a}[y];\nu).
}
$$
Here $\langle\,,\rangle$ is the natural inner product on
$L^2(\varGamma\backslash\r{G})$, and $E_p$ is the 
Eisenstein  series of weight $2p$ (see [2, $(3.19)$]). 
\par
The convergence of the spectral expansion of $\Psi$ is 
in fact absolute and fast, provided $\alpha$ is sufficiently large.
To confirm this, we apply the operator
$\Omega+i\partial^2_\theta$ repeatedly as is done at [2, (5.9)]; and
we get, via the above bound for $a_p(\tau)$,
$$
\langle\Psi(\cdot,\tau),\varphi_p(\cdot,V)\rangle\ll 
(|\tau|+1)^{4\alpha}(|\nu_V|+|p|)^{-{1\over2}\alpha},
$$ 
Then we appeal to the basic bounds, getting
$$
\eqalign{
&\Delta(y,\tau;V)\cr
&\ll (|\tau|+1)^{4\alpha}
(|\nu_V|+1)^{-{1\over4}\alpha}
y^{{1\over2}-\varepsilon}\exp(-\sqrt{y}).
}
$$
The contribution of the continuous spectrum or ra\-ther the
function $\Delta(y,\tau;\nu)$ is to be discussed later.
This yields our assertion. Note that  $\Delta(y,\tau;\cdot)$ 
is regular for $\Re\tau>0$.
\par
\smallskip
In order to pick up a particular Fourier coefficient,
we have used the projection procedure 
as Blomer and Harcos did in [1] to avoid an inner product argument
that was applied in [6]. It is worth pointing it out that in [2,
Sections 3--5] a procedure of  the same kind was developed, employing
the Kirillov model as a main implement to evaluate  
projections  of a certain Poincar\'e series explicitly in terms of its
seed function. Therefore, we surmise that the projections
$\Delta(y,\tau;\cdot)$ could also be handled more precisely with a
modification of the argument of [2]. To this issue we shall
return elsewhere.
\medskip
{\bf 5.} We are now to make the last estimation procedure 
explicit. Thus,  we note that 
$$
\eqalign{
&\langle\Psi(\cdot,\tau),\varphi_p(\cdot,V)\rangle\cr
&={1\over(\overline{\nu_V}^2-{1\over4}-i(2p)^2)^q}
\langle(\Omega+i\partial_\theta^2)^q
\Psi(\cdot,\tau),\varphi_p(\cdot,V)\rangle\cr
&\ll(|\nu_V|+|p|)^{-2q}\!\left\Vert
(\Omega+i\partial_\theta^2)^q
\Psi(\cdot,\tau)\right\Vert
}
$$
for any integer $q\ge0$.
By definition
$$
\Psi(\r{g},\tau)=\sum_{k=1}^\infty\sum_{l=1}^\infty
a_k(\tau)\overline{a_l(\bar\tau)}\varphi_k(\r{g})
\overline{\varphi_l(\r{g})}
$$
with $\varphi_k(\r{g})=\varphi_k(\r{g},A)$,
$\varphi_l(\r{g})=\varphi_l(\r{g},A)$. Since
$$
\Omega=\txt{1\over4}\r{E}\overline{\r{E}}
-\txt{1\over4}\partial^2_\theta
+\txt{1\over2}i\partial_\theta,
$$
with the Maass operator
$
\r{E}=e^{2i\theta}(2iy\partial_x+2y\partial_y-i\partial_\theta),
$
we see that $\Omega \varphi_k\overline{\varphi_l}$ is a linear
combination of $\varphi_{k+j}\overline{\varphi_{l+j}}$, $j=-1,0,1$,
the coefficients of which are polynomials of the second degree
on $k,l$; note that $\nu_A$ is now regarded as a constant. Thus
$$
(\Omega+i\partial_\theta^2)^q
\varphi_k\overline{\varphi_l} =
\sum_{j=-q}^q d_j(k,l;q)\varphi_{k+j}\overline{\varphi_{l+j}},
$$
where $d_j(k,l;q)$ are polynomials of degree $2q$ on $k,l$; and
$$
(\Omega+i\partial_\theta^2)^q\Psi(\r{g},\tau)=
\sum_{k=1}^\infty\sum_{l=1}^\infty
b_{k,l}(\tau;q)\varphi_k(\r{g})
\overline{\varphi_l(\r{g})},
$$
where
$$
\eqalign{
b_{k,l}(\tau;q)&=\sum_{j=-q}^q d_j(k+j,l+j;q)a_{k+j}
(\tau)\overline{a_{l+j}(\bar\tau)}\cr
&\ll (|k|+|l|+1)^{2q}{(|\tau|+1)^{4\alpha}
\over((|k|+1)(|l|+1))^{\alpha}},
}
$$
for each fixed $q$. We then note as before that
$\varphi_k(\r{g})\ll (|k|+1)^2$; and thus,
with $q=\!\left[{1\over3}\alpha\right]$, say, we get the uniform bound
$$
(\Omega+i\partial_\theta^2)^q\Psi(\r{g},\tau)\ll(|\tau|+1)^{4\alpha},
$$
which gives the first inequality in the last section.
\medskip
{\bf 6.}
As to the continuous spectrum, we invoke the
functional equation for $E_p$ ([2, $(3.33)$]), and  have, for
$\Re{\nu}=0$,
$$
\eqalign{
&\langle(\Omega+i\partial_\theta^2)^q\Psi(\cdot,\tau),
E_p(\cdot,\nu)\rangle\cr
&=\gamma_p(\nu)\int_{\varGamma\backslash\r{G}}
(\Omega+i\partial_\theta^2)^q\Psi(\r{g},\tau)E_{-p}(\r{g},\nu)d\r{g},
}
$$
with
$$
\gamma_p(\nu)=\pi^{-2\nu}{\zeta(1+2\nu)\over\zeta(1-2\nu)}
{\Gamma\!\left(\txt{1\over2}+\nu+p\right)\over
\Gamma\!\left(\txt{1\over2}-\nu+p\right)}.
$$
Assuming $\Re\nu>{1\over2}$, we unfold the last integral, and see that
it equals
$R\!\left(\nu+\txt{1\over2}\right)\!Y_p(\nu,\tau;q)/\zeta(2\nu+1)$,
where $R$ is as above, and
$$
\eqalign{
&Y_p(\nu,\tau;q)=
\sum_{l=-\infty}^\infty b_{l+p,l}(\tau;q)\cr
&\cdot\sum_{\delta=\pm}\int_0^\infty\e{A}^\delta\phi_{l+p}(\r{a}[y],
\nu_A)\overline{\e{A}^\delta\phi_l(\r{a}[y],\nu_A)}
y^{\nu-{3\over2}}dy.
}
$$
Again by the basic bounds, we see that $Y_p(\nu,\tau;q)$ is regular
and $\ll(|\tau|+1)^{4\alpha}$ for $\Re\tau>0$ and $\Re\nu>-{1\over2}$. 
Hence, in the same domain,
$$
\Delta(y,\tau;\nu)=
{R\!\left(\nu+\txt{1\over2}\right)\over\zeta(1-2\nu)}
\Lambda(y,\tau;\nu),
$$
with
$$
\eqalign{
\Lambda(y,\tau;\nu)=&\pi^{-2\nu}\sum_{p=-\infty}^\infty
{Y_p(\nu,\tau;q)\over
\left(\nu^2-{1\over4}-i(2p)^2\right)^{q}}\cr
&\cdot{\Gamma\!\left(\txt{1\over2}+\nu+p\right)\over
\Gamma\!\left(\txt{1\over2}-\nu+p\right)}\e{A}\phi_p(\r{a}[y];\nu).
}
$$
One may conclude, via [2, $(4.3)$ and $(4.5)$], that
$$
\eqalign{
&\Lambda(y,\tau;\nu)\cr
&\ll(|\tau|+1)^{4\alpha}
(|\nu|+1)^{-{1\over4}\alpha}
y^{{1\over2}-|\Re\nu|-\varepsilon}\exp(-\sqrt{y})
}
$$
for $\Re\tau>0$ and $\Re\nu>-{1\over2}$.
\medskip
{\bf 7.} We now set $\tau=s$, and observe that
$$
\eqalign{
&D_f(s)={s^{s+2\alpha}\over
\Gamma(s+2\alpha)}\int_0^\infty y^{s-2}\cr
&\quad\cdot\int_0^1
\Psi(\r{n}[x]\r{a}[y],s)\exp(2\pi ifx)dxdy\cr
&=\sum_V f^{{1\over2}-s}\varrho_V(f) \Xi(s,V)\cr
&+\int_{(0)}
{f^{{1\over2}-\nu}\sigma_{2\nu}(f)\over
\zeta(1+2\nu)\zeta(1-2\nu)}
R\!\left(\nu+\txt{1\over2}\right)\Xi(s,\nu){d\nu\over4\pi i},
}
$$
where
$$
\eqalign{
\Xi(s,V)&={s^{s+2\alpha}\over\Gamma(s+2\alpha)}
\int_0^\infty y^{s-2}\Delta(y,s;V)dy,\cr
\Xi(s,\nu)&={s^{s+2\alpha}\over\Gamma(s+2\alpha)}
\int_0^\infty y^{s-2}\Lambda(y,s;\nu)dy.
}
$$
The bound for $\Delta(y,\tau;V)$ implies that
$\Xi(s,V)$ is regular and 
$\ll|s|^{4\alpha+{1\over2}}(|\nu_V|+1)^{-{1\over2}\alpha}$
for $\Re s>{1\over2}$. Similarly, $\Xi(s,\nu)$ is regular
and $\ll|s|^{4\alpha+{1\over2}}(|\nu|+1)^{-{1\over2}\alpha}$
for $\Re s>{1\over2}+|\Re\nu|$ and $\Re\nu>-{1\over2}$.
\smallskip
Therefore we have proved that $D_f(s)$ is indeed regular and of
polynomial growth for $\Re s>{1\over2}$ ([1, Theorem 2]). 
\smallskip
It is to be observed that in order to offset the exponential
growth of the factor $1/\Gamma(s+2\alpha)$ 
we have adopted an old idea
of Ju.V. Linnik which he introduced
in his investigation on approximate functional equations for
Dirichlet $L$-functions. The introduction of the parameter $\tau$
was made with this effect in mind.
\medskip
{\bf 8.} We return to the function $J(u,v;g)$; thus
we impose $\Re(u+v)>\max\{2,1+\eta\}$ initially. In view
of the fast decay of $g^*(s;u,v)$, the last
assertion on $D_f(s)$ yields immediately that 
$$
\eqalign{
&J(u,v;g)=\sum_VL_V\!\left(u+v-\txt{1\over2}\right)\Theta(u,v;V;g)\cr
&+{1\over4\pi i}\int_{(0)}{\zeta(u+v-{1\over2}+\nu)
\zeta(u+v-{1\over2}-\nu)\over\zeta(1+2\nu)\zeta(1-2\nu)}\cr
&\hskip 2cm\cdot R\!\left(\nu+\txt{1\over2}\right)
\Theta(u,v;\nu;g)d\nu, 
}
$$
where
$$
\Theta(u,v;V;g)={1\over2\pi i}\int_{(\eta)}
\Xi(u+v-\xi, V)g^*(\xi;u,v)d\xi,
$$
and the factor $\Xi(u+v-\xi,\nu)$ appears instead in the continuous
spectrum. The $\Theta$ is of fast decay either in $\nu_V$ or
in $\nu$, as $\Xi$ is.
\par
We fix a sufficiently small $\varepsilon>0$; we may move the last 
contour to $(\varepsilon)$, provided
$\Re(u+v)>{2\over3}$, say. Hence, the last expansion of
$J(u,v;g)$ holds under 
$\Re (u+v)>{3\over2}$. This lower bound
is required to
get the factors $L_V\!\left(u+v-{1\over2}\right)$ and
$\zeta\!\left(u+v-{1\over2}\pm\nu\right)$. However,
the former is  entire and of a polynomial order
in $\nu_V$ if $u,v$ are bounded.
Thus the cuspidal part of
$J(u,v;g)$ is regular in a neighbourhood of 
the point $\left({1\over2},{1\over2}\right)$ at which
it takes the value
$$
\sum_V L_V\!\left(\txt{1\over2}\right)\Theta_A(V;g),
$$
with $\Theta_A(V;g)=\Theta\!\left({1\over2},{1\over2};V;g\right)$.
\par
As we are about to deal with the continuous spectrum, 
we should remark that
$\Theta(u,v;\nu;g)$ remains regular in the three complex
variables  and of fast decay in $\nu$, throughout the procedure below,
because of the property of $\Xi(s,\nu)$ mentioned above. Thus, we
restrict $(u,v)$ so that $2>\Re(u+v)>{3\over2}$.
Then, in the last expression for $J(u,v;g)$
one may shift the $\nu$-contour to
$\left({1\over2}+\varepsilon\right)$, encountering the pole at
$u+v-{3\over2}$ with the residue 
$$
-{R(u+v-1)\over\zeta(2(2-u-v))}\Theta\!\left(u,v;u+v-\txt{3\over2}
;g\right)
$$
as well as those of the factor 
$R(\nu+{1\over2})/\zeta(1-2\nu)$; we may assume, without loss of any
generality, that
$u,v$ are such that all the residues in question are finite.
This yields a meromorphic continuation of the continuous spectrum
part,  so that one may move
$(u,v)$ close to $\left({1\over2},{1\over2}\right)$ as far as
$\Re(u+v)>1$ is satisfied; 
this condition is needed to have the last $\Theta$ factor
defined well.  Then, shift the
$\nu$-contour back to the original. All the residual contribution
coming from $R(\nu+{1\over2})/\zeta(1-2\nu)$ cancel out those
arising from the previous shift of the contour. Only the
pole at ${3\over2}-u-v$ contributes newly. The resulting
integral is regular at $\left({1\over2},{1\over2}\right)$; we get the
factor $\Theta_A(\nu;g)=\Theta\!\left({1\over2},
{1\over2};\nu;g\right)$, $\Re\nu=0$.
\medskip
{\bf 9.} Collecting all the above, we obtain
\smallskip
\noindent
{\bf Theorem.} {\it Under the above specifications, 
we have the spectral decomposition
$$
\eqalign{
&I(A;g)=M(A;g)+2\Re\Bigg\{\!\sum_VH_V
\!\left(\txt{1\over2}\right)\Theta_A(V;g)\cr
+&\int_{(0)}\!{\zeta\!\left({1\over2}+\nu\right)\!
\zeta\!\left({1\over2}-\nu\right)\over\zeta(1+2\nu)\zeta(1-2\nu)
}R_A\!\left(\txt{1\over2}+\nu\right)\!\Theta_A(\nu;g)
{d\nu\over4\pi i}\!\Bigg\},
}
$$
where $R_A=R$, and
$M(A;g)$ is the value at $\left({1\over2},{1\over2}\right)$ of
the function
$$
\eqalign{
&{R(u+v)\over\zeta(2(u+v))}
\hat{g}(0)+{1\over\zeta(2(2-u-v))}\cr
&\cdot\Big\{\!R(u+v-1)\Theta\!\left(u,v;u+v-\txt{3\over2};g\right)\cr
&+R(1-u-v)\Theta\!\left(u,v;\txt{3\over2}-u-v;g\right)\cr
&+R(u+v-1)\Theta\!\left(v,u;u+
v-\txt{3\over2};g\right)\cr
&+R(1-u-v)\Theta\!\left(v,u;
\txt{3\over2}-u-v;g\right)\!\Big\}.
}
$$
}
\medskip
\noindent
Albeit the $\Theta$ factors in the last expression is defined so far
only under the condition $2>\Re(u+v)>1$, the expression can
in fact be continued to a neighbourhood of
$\left({1\over2},{1\over2}\right)$, for $I(u,v;g)$ and all other parts
in the spectral expansion of $J(u,v;g)$ and 
$\overline{J(\bar{v},\bar{u};g)}$ are regular
there. We could make the continuation procedure more explicit,
using the property of $\Gamma_p$ mentioned in the first section,
but the above suffices for our present aim.
\smallskip
Our theorem should be compared with [3, Theorem] and [4, Theorem 4.2]
which respectively deal with the mean square of automorphic
$L$-functions associated with discrete series representations and with
that of the product of two values of the Riemann zeta-function,
i.e., the $L$-function associated with the continuous spectrum or
Eisenstein series. Unlike those highly explicit results, admittedly
the above asserts only the existence of a full spectral
decomposition for $I(A;g)$. The construction of the transform
$\Theta_A(\cdot;g)$ has to be made explicit in terms of the weight
function $g$ before one attempts any application; we shall
take this task in our forthcoming works.
Nevertheless, the present work could be a means to broaden the
perspective of the theory of the mean values of the zeta and
$L$-functions that is rendered in [7] via particular examples.
\medskip
{\csc Acknowledgement.} We are greatly indebted to V. Blomer and 
G. Harcos for kindly sending us their important work. 
\bigskip
\centerline{\bf References}
\medskip
\item{[1]} V. Blomer and G. Harcos.  The spectral decomposition
of shifted convolution sums. Preprint, March 8, 2007.
\item{[2]} R.W. Bruggeman and Y. Motohashi. A new approach to the
spectral theory of the fourth moment of the Riemann zeta-function. 
J. reine angew.\ Math., {\bf 579} (2005), 75--114. 
\item{[3]} Y. Motohashi. The mean square of Hecke $L$-series attached
to  holomorphic cusp-forms. Ko\-kyuroku, RIMS Kyoto Univ., {\bf 886}
(1994), 214--227.
\item{[4]} Y. Motohashi. Spectral Theory of the Riemann Zeta-Function.
Cambridge Univ.\ Press, Cambridge 1997.
\item{[5]} Y. Motohashi. A note on the mean value of the zeta and
$L$-functions.\ XII. Proc.\ Japan Acad., {\bf 78A} (2002), 36--41.
\item{[6]} Y. Motohashi. A note on the mean value of the zeta and
$L$-functions.\ XIV.  Proc.\ Japan Acad., {\bf 80A} (2004), 28--33.
\item{[7]} Y. Motohashi. Mean values of zeta-functions via
representation theory. Proc.\ Symp.\ Pure Math., AMS, 
{\bf 75} (2006), 257--279.
\bye